\documentclass[12pt]{article}
\usepackage{amssymb, url, amsmath, amsthm, curves}

\def\Pr{\operatorname{Pr}}

\newtheorem{theorem}{\noindent Theorem}

\newtheorem{definition}{\noindent Definition}
\newtheorem{corollary}{\noindent Corollary}

\newtheorem{statement}{\noindent Proposition}

\urldef\vershik\url{vershik@pdmi.ras.ru}
\author{A.~M.~Vershik\thanks{%
St.~Petersburg Department of Steklov Institute of Mathematics.
E-mail: \vershik.
Partially supported by RFBR, grant 02-01-00093, and  INTAS, grant
03-51-5018.
}}
\title{Polymorphisms, Markov processes, quasi-similarity.}

\date{November 26, 2004}

\begin{document}

\maketitle

\flushright{\it To the centenary of my mother H.~J.~Lusternik.}
\bigskip

\abstract{In this paper we develop the theory of {\it polymorphisms} of measure
   spaces, which is a generalization of the theory of measure-preserving
   transformations. We describe the main notions
   and discuss relations to the theory of markov processes,
   operator theory, ergodic theory, etc. We formulate the important
   notion of quasi-similarity and consider quasi-similarity
   between polymorphisms and automorphisms.

   The question is as follows: is it possible to have a quasi-similarity
   of a measure-preserving automorphism  $T$ and a polymorphism $\Pi$ (which
   is not an automorphism)? In less definite terms: what kind of
   equivalence can exist between deterministic and
   random (markov) dynamical systems?
   We give the answer: each nonmixing prime polymorphism is quasisimilar to
   the automorphism with positive entropy and  for each $K$-automorphism $T$ is
   quasisimilar to a polymorphism $\Pi$ which is a special random
   perturbation of the automorphism $T$.}
\newpage
 \tableofcontents

\section{Introduction}
 \subsection{The theory of polymorphisms as a generalization of
 ergodic theory}

   The simplest example of a polymorphism with invariant measure
   that is not a measure-preserving auto- or endomorphism
   is the ``map'' from the unit circle
   $S^1$ with Lebesgue measure to itself
   whose graph is the following cycle (curve) on the
   2-torus $S^1\times S^1$ (also with Lebesgue measure):
   $$\{(u,v):u^2=v^3,\,|u|=|v|=1\}.$$
   Another example is a ``random automorphism,'' or the convex combination
   $p_1 T_1 + \dots + p_n T_n$, where $T_i$, $i=1,\dots, n$,
   are automorphisms of a measure space, $p_i >0$, $\sum p_i=1$.

   In these examples,
   each point of the space under consideration has
   many images (a random image) and a random preimage.
   A general polymorphism with invariant measure from one measure space
   to another is, by definition, a {\it measure in the product of these
   spaces with given marginal projections} (see Definition 1). It was first defined
   in the paper \cite{Vpol}. The dynamics of such ``maps'' is an interesting
   open area.

   In the first part of the paper we give a brief survey of definitions and
   important properties of the notion of polymorphism. Briefly speaking, we
   {\it generalize the foundations of ergodic theory to a natural measure-theoretic
   version of the dynamics of multivalued maps}. We called the main objects of
   our theory  {\it polymorphisms}; a polymorphism can be regarded simply as a
   {\it Markov map with a fixed invariant measure},
   or the two-dimensional distribution of a stationary Markov process, or
   a bistochastic measure, or a joining of measure spaces.

   Parallel notions to the notion of polymorphism in other parts of mathematics
   are: {\it correspondence} in algebra and algebraic geometry;
   {\it bifibration} in differential geometry, {\it Markov map} in probability theory,
   {\it Young measure} in optimal control, etc. The notion of
   polymorphism (also with a quasi-invariant measure) generalizes all
   such examples, see Section 2.3 (``Why polymorphisms?''), but we are mainly
   interested in the geometry and dynamics of polymorphisms in the
   framework of {\it measure theory}.

   From the point of view of dynamics and physics, the notion of polymorphism
   corresponds to the {\it ``coarse-graining'' approach} to dynamics:
   instead of one-to-one maps we are allowed to consider maps that send a point
   to a measure, or elements of a partition to a family of elements of the same
   partition (a set of ``grains'' to itself, see Section 2).
   This opens new possibilities, which are forbidden in the classical theory: for
   example, we can define the notion of the quotient of an automorphism by a
   partition {\it that is not invariant under the automorphism}; this quotient is a
   polymorphism, but not an automorphism. Another important direction is approximation
   of automorphisms with polymorphisms, especially with finite
   polymorphisms; this approach is an alternative to Rokhlin approximations and
   is perhaps more effective for automorphisms with positive entropy.

\bigskip
    The set of polymorphisms of a given measure space has a rich structure:
    it is a convex weakly compact topological semigroup whose invertible
    elements are measure-preserving automorphisms.
    The functional analog or operator formalism of our theory is the theory of
    Markov operators in the Hilbert space $L^2_{\mu}(X)$. {\it A Markov operator is
    a positive contraction that preserves constants};  the positivity means preserving
    the cone of nonnegative functions in the space $L^2_{\mu}(X)$.
    This immediately leads us to the necessity of generalizing the theory of contractions
    and non-self-adjoint operators in Hilbert spaces to Markov operators.

    The spectral theory of concrete Markov operators had been studied for a long
    time; but the theory of contractions in Hilbert spaces apparently did
    not up to now attract the attention of specialists in dynamical systems;
    it gives a further development of spectral theory and a new type of
    questions in dynamics, which is essentially important for polymorphisms.
    We choose one of such questions that is interesting in itself
    and describe it briefly
    in the next section.

    \subsection{The problem of quasi-similarity of automorphisms and polymorphisms,
    and paradoxical Markov processes}

    Recall (see \cite{FN}) that a bounded operator $T$ in a Hilbert space $H_1$
    is called a {\it quasi-image} ("quasiaffinitet" in \cite{FN}), of a bounded operator
    $S$ in a Hilbert space $H_2$
    if there exists a continuous linear operator $L$ from  $H_1$ to $H_2$ that
    may have no bounded inverse but has a \textbf{dense image}
    (so that the inverse operator is defined
    on a dense set) such that $LT =SL$. Two bounded operators are called
    {\it quasisimilar} if each of them is a quasi-image of the other.
    Quasi-similarity is an equivalence in the space of bounded operators, which
    is weaker than similarity (or unitary equivalence) of operators;
    it may happen that there is an equivalence between a unitary operator and a proper
    contraction operator that is totally nonunitary (see \cite{FN}). This
    equivalence does not in general preserve spectra of operators.
    Such examples are important in functional analysis and scattering theory.

    Quasisimilarity for Markov operators and the parallel geometric notion
    of quasi-similarity for measure-preserving transformations and especially for
    polymorphisms seems to have never been considered systematically.
    A polymorphism $\Pi_1$ is a quasiimage of a polymorphism $\Pi_2$
    (in particular, one or both of them can be automorphisms)
    if there exists a dense
    (in the sense of Section 2) polymorphism $\Lambda$ such that
    $\Lambda\cdot \Pi_1 =\Pi_2\cdot \Lambda$. Equivalently this means that
    there exist joining between $\Pi_1$ and $\Pi_2$.
    {\it Two polymorphisms
    are quasisimilar if each of them is a quasi-image of the other.}
    The question is as follows: can such a $\Lambda$ exist if $\Pi_1$ is an automorphism
    and $\Pi_2$ is a proper polymorphism? The general problem is
    \textbf{to describe all quasi-similar pairs ``automorphism $\leftrightarrow$
    po\-ly\-mor\-phism.''}

\medskip

    In other words, the problem is to describe pairs ``deterministic
    transformation $\leftrightarrow$ random transformation'' that can be quasi-similar
    in the above sense. It may happen that this question can be related
    to the long discussion among physicists (see, e.g., \cite{Pr, Mis, Gust}) on
    possible equivalence between deterministic and random systems.

    In a nontrivial case, a polymorphism must be ``prime,'' i.e., have no
    factor endomorphisms, and also non-mixing in the sense of the theory of Markov
    processes; otherwise the problem is not interesting. Even the existence
    of such polymorphisms is not obvious.
    The first example was given in probabilistic terms by M.~Rosenblatt
    (see \cite[Ch.~4.4]{Ros}). Then a smooth example was suggested
    in \cite{Vgrub}; it was a random perturbation of a hyperbolic
    automorphism of the torus. In this paper we give a formulation
    of the complete solution of the problem.

    Here we present one result in this direction and leave the general case
    for another publication (see Section 5).

    {\it 1.There is the automorphism $T$ with positive entropy which a canonical quasiimage
    of  non-mixing prime polymorphisms $\Pi$; if conjugate polymorphism $\Pi^*$ is also prime
    and nonmixing then polymorphism $\Pi$ is quasisimilar to $T$.

    2.Assume that for $K$-automorphism $T$: there exists a finite
    or countable $K$-generator such that in the symbolic realization of $T$
    with this generator, the homoclinic equivalence relation is ergodic,
    then there exists a polymorphism $\Pi$ that
    is quasi-similar to $T$, more exactly, the following weak limits exist and
    define two intertwining polymorphisms
    $$\Lambda_1=\lim_{n\to \infty}\Pi^n T^{-n},\qquad \Lambda_2=
    \lim_{n\to \infty}T^n(\Pi^*)^n$$
    that realize the quasi-similarity:
      $$\Pi \cdot \Lambda_1 = \Lambda_1 \cdot T, \qquad
      \Lambda_2 \cdot \Pi =\Pi \cdot \Lambda_2.$$
     The polymorphism $\Pi$ is a special random perturbation of the
     $K$-automorphism $T$.}

     The question of whether this $K$-automorphism $T$ is unique
     leads to a very interesting problem for automorphisms and especially
     $K$-automorphisms: \textbf{\it do there exist two automorphisms that
     are not isomorphic but are quasi-similar?} I do not know the
     answer.

     In order to explain our method, we must say several words
     about stationary Markov processes with paradoxical property
     that appeared in this problem.

     It is well known that a mixing Markov chain is regular in the
     sense of Kolmogorov (or pure nondeterministic): the
     $\sigma$-field of the infinite past is trivial. On the other hand,
     mixing is equivalent to irreducibility, i.e., the absence of
     nontrivial partitions of the state space into subclasses,
     or the absence of deterministic factors of a Markov process
     (or primality for polymorphims).\footnote{Perhaps, this
     fact was first proved in the paper \cite{Vin} by a
     student of Kolmogorov.}

     Contrary to this, the absence of nontrivial (measurable)
     partitions into subclasses, or the absence of nontrivial deterministic
     factor processes, which we called primality for polymorphisms,
     does not imply mixing and regularity for general Markov processes
     --- there exist nonmixing Markov processes that have no deterministic
     factors; we called such paradoxical Markov processes
      {\it quasi-deterministic} Markov processes.

     \smallskip
     The reason of this difference between processes with discrete and continuous
     state spaces is rather deep and relates to the
     theory of measure-theoretic equivalence relations.
     A {\it measurable} partition of the state space of a Markov chain
     into ``subclasses'' allows us to decompose the Markov chain into irreducible
     mixing (nonhomogeneous in time) chains. But in the general case it may happen
     that there exists a {\bf nonmeasurable} partition of the state space,
     or an ergodic equivalence relation that is invariant under the
     polymorphism, so that there is no regularity, but at the same time there
     are no deterministic factors. This effect underlines
     quasi-similarity. As we will see, such processes have a
     hidden determinism: you cannot predict
     the value of the process at time zero  with probability one
     if you know the infinite
     past, but you can define a conditional probability on this
     state space, and this conditional probability is
     different for different points of the tail space
     (``entrance boundary''). The action of the shift on the tail
     $\sigma$-field gives an action on the set of these conditional
     measures.

     The role of the $K$-automorphism mentioned in the theorem is
     played by the ``tail shift'' --- the restriction of the shift
     onto the tail $\sigma$-field;
     the Markov generator is not a $K$-generator for the Markov shift; nevertheless,
     this is a $K$-shift. In order to prove the $K$-property, we must
     change the generator. All details will be published in a
     separate paper.

\smallskip
     In the second section we give the first definitions, examples, links, etc.
     In the third section we introduce the operator
     formalism and operator version of quasi-similarity. Section 4
     is devoted to the corresponding Markov processes,
     tail $\sigma$-fields, residual automorphisms. The main results
     are formulated in Section 5.

     Note that we consider the case when the time is $\Bbb Z$, but there
     are no serious obstacles to extending the results to the continuous
     time $\mathbb R$. We will return to this topic in more detail elsewhere.

\section{Definitions and properties of polymorphisms}

  We will briefly define the main notions we need.
  Some details can be found in \cite{Vpol}.

    The notion of polymorphism is a measure-theoretic analog
    of what people called a multivalued map. In the framework of measure
    theory, the value of a ``multivalued map'' at a point is not a subset
    of the target space, but a measure on this space. In this sense, a
    polymorphism is a measure-theoretic analog of a Markov map; in the
    subsequent sections we will discuss the relation to the
    theory of Markov processes in detail.

    Objects similar to polymorphisms have many names
    in various theories (see the introduction).
    Our considerations are directed towards {\it dynamics, probability,
    and ergodic theory.}

\subsection{First definitions}

   Let $(X,\mu)$ be a Lebesgue space with continuous measure $\mu$
  (i.e., a measure space isomorphic to the unit interval with the Lebesgue measure).

 \begin{definition}
 {A polymorphism $\Pi$ of the Lebesgue space $(X,\mu)$ to itself
    with invariant measure $\mu$ is a diagram consisting of an ordered triple
    of Lebesgue spaces:
    $$ (X,\mu)  \stackrel{\pi_1}\longleftarrow (X \times X, \nu)
    \stackrel{\pi_2}\longrightarrow (X, \mu),$$
    where $\pi_1$ and $\pi_2$ stand for the projections to the first
    and second component of the product space $(X \times X, \nu)$, and
    the measure $\nu$, which is defined on the $\sigma$-field generated by the
    product
    of the $\sigma$-fields of $\bmod0$ classes of measurable
    sets in $X$, is such that $\pi_i\nu=\mu$, $i=1,2$.}

    The measure $\nu$ is called the
    {\it bistochastic measure of the polymorphism $\Pi$}.
\end{definition}

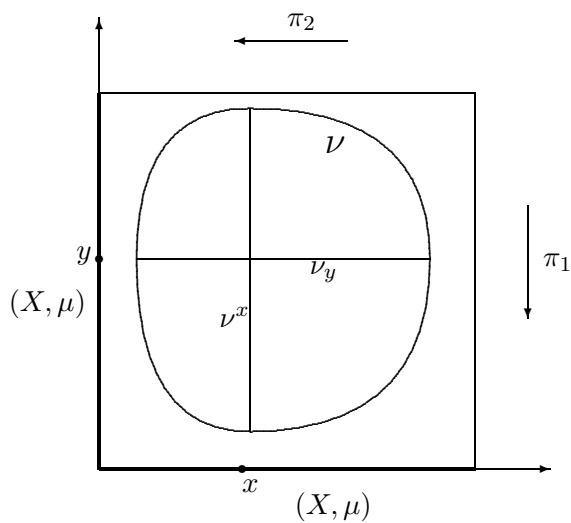
\setlength{\unitlength}{1mm}
\begin{center}
\begin{picture}(60,60)
\put(0,0){\vector(1,0){60}} \put(0,0){\vector(0,1){60}} \put(57,35){\vector(0,-1){15}}
\put(59,27){$\pi_1$} \put(33,57){\vector(-1,0){15}} \put(25,59){$\pi_2$} \put(20,5){\line(0,1){43}}
\put(19,-3){$x$} \put(19,0){\circle*{1}} \put(50,0){\line(0,1){50}} \put(5,28){\line(1,0){39}}
\put(-3,28){$y$} \put(0,28){\circle*{1}} \put(0,50){\line(1,0){50}} \closecurve(20,5,  44, 28,
20,48,  5, 28) \put(16,19){$\nu^x$} \put(28,26){$\nu_y$} \linethickness{.4mm}
\put(0,0){\line(0,1){50}} \put(0,0){\line(1,0){50}} \put(30,42){\large $\nu$}
\put(26,-6){$(X,\mu)$} \put(-12,21){$(X,\mu)$}
\end{picture}
\bigskip\bigskip\bigskip
\end{center}

    A polymorphism $\Pi^*$ is called {\it conjugate} to the polymorphism
    $\Pi$ if its diagram is obtained from the diagram of  $\Pi$
    by reflecting with respect to the central term.

    Consider the ``vertical'' partition $\xi_1$ and
    the ``horizontal'' partition $\xi_2$ of the space $(X\times X, \nu)$
    into the preimages of points under the projections $\pi_1$ and $\pi_2$,
    respectively.
    In terms of bistochastic measures, the
    {\it value of a polymorphism at a point $x\in X$}
    is a conditional measure. More precisely, we have the following definition.

\begin{definition}
   {In the above notation, the value of the polymorphism
   $\Pi:X \to X $ at a point $x_1 \in X$
   is, by definition,  the conditional measure $\nu^{x_1}$
   of $\nu$ on the set $\{(x_1,\cdot)\}$ with respect to the vertical
partition $\xi_1$
   (the transition probability); similarly,
   the value of the conjugate polymorphism $\Pi^*$ at a point $x_2\in X$
   is the conditional measure $\nu_{x_2}$
   of $\nu$ on the set $\{(\cdot,x_2)\}$ with respect to the horizontal
   partition $\xi_2$ (the cotransition probability).
      These conditional measures
   are well-defined on sets of full measure.}
  \end{definition}

\smallskip\noindent{\bf Remark.}
It is very important that the set of conditional measures $\{\nu^x, \,x \in
X\}$ itself does not determine the polymorphism --- we need to know also the
measure $\mu$ on $X$. Recall that an ordinary Markov map
is determined by the list of transition probabilities.
\smallskip

    With obvious modifications, we can define more general notions:

    1) a polymorphism of one measure space $(X_1,\mu_1)$
    to another measure space $(X_2,\mu_2)$:
    $$ (X_1,\mu_1)  \longleftarrow (X_1 \times X_2, \nu)
    \longrightarrow (X_2, \mu_2),$$ where the measure $\nu$
    have marginal projections $\mu_1$ and $\mu_2$;

    2) a polymorphism with {\it quasi-invariant measure}; in this case,
    the projections $\pi_1\nu$ and $\pi_2\nu$
    of the measure $\nu$ are equivalent (not necessarily equal)
    to the measures $\mu_1$ and $\mu_2$, respectively.

    {\it For the most part, we will consider polymorphisms
    of a space with continuous measure to itself  with a finite invariant measure.}

    All notions should be understood
    up to set of zero measure ($\bmod0$).
    In fact, our objects and morphisms are classes of
    coinciding $\bmod 0$ objects and morphisms, but we will not repeat
    the corresponding routine
    comments when this does not cause any problem.

    For simplicity, we assume that
    the $\sigma$-field $\emph{A}$ of the Lebesgue space $(X,\mu)$
    has a countable basis $\emph{B}$.
    This means that on $X$
    we have the standard Borel structure;
    consequently, on the space $X\times X$
    we have the countable basis $\emph{B}\times \emph{B}$,
    the standard Borel structure, and the
    $\sigma$-field generated by this basis. Thus all bistochastic measures
    $\nu$ corresponding to polymorphisms of
    the space $(X,\mu)$ to itself
    will be defined on this $\sigma$-field.\footnote{We omit
    the discussion of the nontrivial question concerning
    the independence of such a definition on the choice of the basis
    $\emph{B}$.}  Now the set of all bistochastic measures becomes an
    affine compact space equipped with the topology of weak convergence
    on the basis.

\subsection{Further definitions and properties}

The following proposition-definition describes structures on
polymorphisms.

\begin{statement} {The set of polymorphisms
    (bistochastic measures) $\cal P$
    is a topological semigroup with the
    following natural product. Let $\Pi_1,\Pi_2$ be two polymorphisms
    with bistochastic measures $\nu_1,\nu_2$; then the
    product $\Pi_1\Pi_2$  has bistochastic measure $\nu$ defined by
    $$\nu^x(A)=\int \nu_1^y(A)\,  d\nu_2^x(y).$$

    The semigroup $\cal P$ has the {\it zero element}: this is
    the polymorphism $\Theta$ with bistochastic measure $v=\mu \times \mu$; obviously,
    $\Theta \cdot \Pi = \Pi \cdot \Theta = \Theta$
    for every $\Pi$; we call $\Theta$ the zero polymorphism.

    The set of polymorphisms $\cal P$ has the structure of a semigroup
    with involution $\Pi \to \Pi^*$, which was defined above.

    The subgroup of invertible elements of the semigroup $\cal P$
    is the group of measure-preserving transformations;
    the semigroup of measure-pre\-ser\-ving
    endomorpisms is a subsemigroup of $\cal P$;
    the bistochastic measure corresponding to an endomorphism $T$
    is the measure on the set $\{(x,Tx)\}_{x \in X} \subset (X \times X)$ that is
    the natural image of the
    measure $\mu$ under the map $x \to (x,Tx)$.}
\end{statement}

All assertions of the proposition are obvious. We explain only the last one.
Assume that $T$ is an endomorphism of a Lebesgue space $(X,\mu)$
    with invariant measure. Consider the graph of $T$,
    i.e., the set $\{(x,Tx)\}_{x \in X} \subset X \times X$, and the measure
    $\mu$ on this graph (more rigorously, we identify a point $x \in X$
    with the point $(x,Tx)$, so that the measure $\nu$ can be
    regarded as the image
    of the measure $\mu$ on the graph of $T$). Since $T$
    is a measure-preserving map, it follows that $\nu$ is a bistochastic measure,
    and we can identify the endomorphism $T$ with the
    corresponding polymorphism.
    Thus we embed the semigroup of endomorphisms and, in particular,
    the group of automorphisms into the semigroup of polymorphisms.

\smallskip
    Let us define the notions of factor polymorphism, ergodicity, mixing,
    etc., and compare them with the same notions for endomorphisms.

\medskip

  1. A measurable partition $\xi$ is called {\it invariant} under a
  polymorphism $\Pi$ if for almost all elements $C \in \xi$ there exists
  another element $D \in \xi$ such that for almost all (with respect to the
  conditional
  measure on $C$) points $x \in C$, we have $\mu^x(D)=1$, where $\mu^x$ is the $\Pi$-image of $x$.
  In other words, the factor polymorphism $\Pi_{\xi}$ of $\Pi$
  by an invariant partition $\xi$ is an endomorphism of the space $(X_{\xi},\mu_{\xi})$.

  In particular, if for almost all elements $C \in \xi$ of
  $\xi$ we have $\mu^x(C)=1$ for almost all $x\in C$, then the
  partition $\xi$ is called a {\it fixed} partition for $\Pi$
  and the corresponding
  factor polymorphism is the identity map on $X_\xi$.

Both definitions restricted to endomorphisms give the corresponding
  notions (of invariant and fixed
  partitions) of ergodic theory. Any polymorphism has the maximal fixed partition
  and can be decomposed into the direct integral  of ergodic components
  over this partition.

\smallskip
  2. The {\it ergodicity} of a polymorphism $\Pi$
  means the absence of fixed measurable
  partitions, in other words, the absence of identical factors.

\smallskip
  3. Given a polymorphism $\Pi$ of
  a space $(X,\mu)$ with bistochastic measure $\nu$,
    the {\it factor polymorphism} of $\Pi$
  by a measurable partition $\xi$
  is the polymorphism of the space $(X/{\xi},\mu/{\xi})$ to itself
  with bistochastic measure $\nu/{(\xi \times \xi)}$.
   That is,  we have the diagram
  $$ (X_\xi,\mu_\xi)  \longleftarrow (X_\xi \times X_\xi, \nu_{\xi \times \xi})
    \longrightarrow (X_\xi, \mu_\xi).$$

  Thus the factor polymorphism of any polymorphism by any measurable
  partition does exist, in particular, the {\it factor polymorphism of
  any automorphism by any (not necessarily invariant) partition
  always exists}.

\smallskip
 4. A polymorphism is called {\it prime}\footnote{In \cite{Vgrub},
   this notion was called ``exactness.''} if it has no nontrivial
   invariant partitions.\footnote{There is a small difference
   between this notion and
   the notion of a prime automorphism, which is, by definition, an
   automorphism that has no
   invariant partitions except the trivial partition $\nu$ and the
   partition into separate points
   $\varepsilon$; thus a prime automorphism is not a prime polymorphism in
   our sense, because $\varepsilon$ is an invariant partition for
   it.} Prime nonmixing polymorphsms are the main objects
   of the second part of this paper. In a sense, primality is an analog of Rokhlin's notion of
   exactness for endomorphisms.

\smallskip
  5. A polymorphism $\Pi$ is called {\it mixing} if the sequence of its powers
   tends to the zero polymorphism (see above) in the weak topology:
   \mbox{w-$\lim_{n \to \infty}\Pi^n = \Theta$}.  Note that
   it may happen that a polymorphism
   is mixing while its conjugate is not.

   This notion of mixing has nothing to do with the
   notion of mixing in ergodic theory:
   for example, no automorphism is mixing in our sense. We use this term,
   because it is equivalent to the traditional notion of mixing in the sense of the theory of
   Markov processes --- see the next section.

\smallskip
   6. We say that a polymorphism $\Pi$ is
   {\it injective} if the partition into the preimages
   of points under the map $x \to \Pi(x)$ coincides with $\varepsilon$ (= partition into
   separate points $\bmod 0$). We will deal with injective
   polymorphisms in Section 5.
   If the measures $\Pi(x)$ are discrete
   (finitely supported)
   for almost all points $x$, we say that the polymorphism
   $\Pi$ is of {\it discrete rank}
   (respectively, of {\it finite rank}).
   If the bistochastic measure of a polymorphism $\Pi$
   of a space $(X,\mu)$ is
   absolutely continuous with respect to the product measure
   $\mu\times\mu$,
   we say that
   $\Pi$ is {\it absolutely continuous}.

\smallskip
  7. A polymorphism $\Pi$  called {\it dense} if there is no nonzero measurable function $f$
  that has zero mean with respect to $\mu$-almost all (in $x$) conditional
  measures $\Pi(x)=\nu^x$, that is,
  if $\int f(y)d\nu^x(y)=0$ for $\mu$-almost all $x$ implies $f\equiv 0$.
  Below we will give an equivalent definition.

\medskip
   We will not continue the list of definitions and restrict ourselves only with
   notions we need in this paper. For example, we do not consider the entropy of
   polymorphisms, spectral properties, etc.

\subsection{Why polymorphisms?}

      The notion of polymorphism widely extends the theory of
   transformations with invariant measure. We briefly illustrate
   some advantages
   and aspects of this notion.

\medskip
   \textbf{1. A polymorphism of a finite measure space} $X$ with the uniform measure
   $\mu(\cdot)=1/\#X$ to itself is a bistochastic matrix
   of order $m=\#X$. Such matrices form the semigroup
   $V_m=\{(a_{i,j})_{i,j=1}^m: \sum_j
   a_{i,j}=\sum_i a_{i,j}=1/m,\, a_{i,j} \geq 0\}$.\footnote{It is convenient
   to assume that
   the rows of matrices sum to $1/m$ and not to $1$ as usual.}
   Thus the factor polymorphism of any polymorphism by a finite partition
   with $n$ parts can be identified with a bistochastic matrix of order $n$.

  This leads to the following easy proposition.

\begin{statement}
   {Every polymorphism of a continuous measure space can be represented
   as the inverse limit of a sequence of polymorphisms of finite
   spaces with uniform measures, or, in other words, of bistochastic
   matrices, for example, of orders $2^n$. The semigroup ${\cal P}$ of all polymorphisms
   of a continuous measure space with the weak topology is the inverse limit
   of a sequence of semigroups of bistochastic matrices: ${\cal P}= \projlim_n V_{2^n}$.}
\end{statement}

   This gives an alternative approach to approximation
   in the ergodic theory of automorphisms. For some reasons, approximation
   of automorphisms with positive entropy by bistochastic matrices
   (= ``periodic'' polymorphisms) is more natural than approximation by
   periodic automorphisms. Following the physical terminology, we can call
   this type of approximation the {\it coarse-graining approximation
   of dynamical systems} (see \cite{Vpol}).
   I think that it is a fruitful method of studying $K$-automorphisms.

\smallskip

  \textbf{2. Polymorphisms allow us to extend  ordinary notions in a natural way}.
   For example, the conjugate to an endomorphism does not exist in the ordinary
   sense, but does exist as a polymorphism; we called such a polymorphism an
   {\it exomorphism}: in this case, a point has several
   images but one preimage.

\smallskip
   As we have seen, the notion of polymorphism allows us to consider quotients
   of arbitrary automorphisms by arbitrary measurable partitions.

   We may say that the theory of polymorphisms is the envelope
   of the theory of endomorphisms with respect to extending
   the notion of factorization from invariant partitions to
   arbitrary ones.

   Moreover, any two partitions, say $\zeta$ and $\eta$, produce
   a polymorphism with invariant  measure:
     $$(X_{\zeta},\mu_{\zeta}) \longleftarrow (X_{\zeta \vee \eta},
     \nu_{\zeta \vee \eta}) \longrightarrow
     (X_{\eta},\mu_{\eta}),$$
     where all three spaces are the quotients of the space $(X,\mu)$
    by the corresponding partitions. We will use this remark later.

    In the special case when we have one partition $\xi$ (and the second one is the
    partition into separate points) we have the
    ``tautological'' polymorphism from $X/\xi$
    to $X$, which associates with an element $C \in X/\xi$
     the conditional measure on $C \in X$.

\smallskip
   \textbf{3. Polymorphisms and correspondences}.
   The simplest classical source of polymorphisms with finitely many
   images and preimages is a {\it correspondence} in the sense of algebraic geometry;
   the following scheme gives the simplest example.
   Consider the 2-torus $\{{\Bbb T}^2=(u,v): u,v \in \Bbb C,\, |u|=|v|=1\}$
   and the curve $$u^n=v^m,\quad n,m >1,$$
   equipped with the Lebesgue measure. The conditional measures
   (transition and cotransition) are the uniform measures.
   We obtain a polymorphism of the unit circle with the Lebesgue measure to itself.
   The dynamics of such polymorphisms is very interesting and still
   poorly studied. This example can also be regarded as a bifibration
   over the circle.

\smallskip
 \textbf{4. Polymorphisms and Markov processes}.
  Below we will describe the link to the theory of Markov
  processes: a polymorphism is the
  two-dimensional distribution
  of a stationary Markov process.
  The theory of polymorphisms leads to a new kind
  of examples and problems (see the second part of the paper)
  and help to understand the structure of general Markov
  processes.

\smallskip
  \textbf{5. Random walks on automorphisms as polymorphisms.}
  Another typical example of polymorphisms came from the theory
  of random walks: in this case, a polymorphism is a
  (finite or infinite) convex combination
  of deterministic transformations, for example,
  shifts on some group of measure-preserving transformations;
  the coefficients of this convex combination may depend on
  points: assume that $\{T_{\alpha},\, \alpha \in A\}$ is a family
  of transformations with quasi-invariant measure $\mu$; then
  $\Pi(x)=\mu^x$, where $\mu^x$ is a measure on the set
  $\{T_{\alpha}x: \alpha \in A\}$, or, better to say, $\mu^x$
  is a measure on the set of parameters $A$ that depends on $x$.
  For $\Pi$ to be a polymorphism with invariant measure,
  these measures must satisfy some conditions.

\smallskip
   {\textbf6.} The theory of \textbf{polymorphisms} is closely related to the theory
   of \textbf{joinings}, which are nothing more than polymorphisms with
   additional symmetries (for example, commuting with the automorphism
   $T\times T$ of the space $(X\times X,\nu)$ in the above notation).
    We can also say that a polymorphism with invariant measure is a joining
   of identical maps. It is more important that the quasi-similarity of
   two automorphisms can be also formulated as a problem on
   joinings of special type.

\smallskip
   \textbf{7. Orbit partition of a polymorphism}.
    The {\it trajectory partition}, or {\it orbit partition}
    of a polymorphism $\Pi$ of a space
    $(X,\mu)$
    is defined as follows: two points $x,y$  belong
    to same orbit if and only if there exist
    positive integers $n,m$ such that the measures
    $\Pi^n(x)$ and $\Pi^m(y)$ are not mutually singular as measures on $X$.
    Denote by $o(x)$ the orbit of a point $x$ under the polymorphism $\Pi$.
    If the polymorphism is of discrete rank (see the definition above),
    then the orbit partition has countable fibers. In this case,
    we obtain a new wide class of nonmeasurable partitions, or ergodic
    equivalence relations;
    a very intriguing question is to
    find a criterion of hyperfiniteness
    (tameness) of these partitions or to study their properties
    in terms of polymorphisms.

\smallskip
   \textbf{8. Polymorphisms and groupoids.} This is a very important
   link. For simplicity, assume that an ergodic polymorphism $\Pi$ is of discrete
   rank. Then its orbit partition defines an ergodic equivalence
   relation and a measurable groupoid (see \cite{Ren}).
   We call $\Pi$ {\it complete} if the measure
   $\Pi(x)$ is strictly positive on the orbit $o(x)$ for almost all $x$;
   in this case, the $\Pi$-image of $x$ is a measure on the whole orbit of $x$, or, in other words, the
   bistochastic measure of $\Pi$ is positive on the groupoid.

\subsection{Quasi-similarity of polymorphisms and automorphims}

   A classification of polymorphisms (or Markov operators, see below)
   can be defined in many ways. One of them is the classification up to
   conjugation with respect to
   a measure-preserving automorphism (see \cite{Vpol} for discussion).
   In this paper we will consider
   the classification up to quasi-similarity.

  \begin{definition}
{A polymorphism (in particular, auto- or endomorphism) $\Pi_1$ is
a quasi-image of a polymorphism $\Pi_2$ if there exists a dense
polymorphism $\Gamma$ such that
     $$\Gamma \cdot \Pi_1= \Pi_2 \cdot \Gamma.$$

 We say that two polymorphisms are quasi-similar if each of them is
 a quasi-image of
 the other.}\footnote{The density of $\Gamma$ in this definition
 is a very important condition; without it, the equivalence is trivial.}

  \end{definition}

  \textbf{Question}. To describe the notion of quasi-similarity for
  measure-preserving
  auto- and endomorphisms: does it coincide with the notion of isomorphism? It is especially
  important to know the answer for $K$-automorphisms.

\smallskip
  But we will study the special case when $T$ is a measure-preserving automorphism
  of a Lebesgue space $(X,\mu)$ and $\Pi$ is a polymorphism with invariant
  measure of the same space.

\smallskip
   \textbf{Problem}: {\it When does exist a dense polymorphism $\Gamma$ of $(X,\mu)$
   such that
              $$\Gamma \cdot \Pi = T \cdot \Gamma\quad?$$}
   A similar problem: {\it when does exist a dense polymorphism $\Lambda$ such
   that
              $$\Pi \cdot \Lambda = \Lambda \cdot T\quad?$$
  Or when the automorphism $T$ and the polymorphism $\Pi$ are quasi-similar?}
\smallskip

Looking ahead and using the notions that will be introduced later,
   we can say that, in order to avoid trivial cases
   (when both $T$ and $\Pi$ are automorphisms),
   we should suppose that $\Pi$ is {\it prime} (= has no nontrivial factor
   endomorphisms, or has no nontrivial invariant partitions).
   Furthermore, a mixing polymorphism cannot be quasi-similar to any measure-preserving
   transformation of a continuous measure, thus we may assume
   without lost of generality
   that $\Pi$  (or $\Pi^*$) is nonmixing;
   this means that $\Pi^n \nrightarrow \Theta$
   (respectively, $\Pi^{*n}\nrightarrow \Theta$) in the weak topology
   as $n \to \infty$,
   where $\Theta$ is the zero polymorphism.

    \section{The operator formalism and Markov operators}

   In this section we consider two alternative languages for the theory
   of polymorphisms: the first one is the {\it operator formalism}
   in the space of measurable square
   integrable functions $L^2_{\mu}(X)$,
   the language of so-called {\it Markov operators},
   and the second one is the language of {\it stationary Markov processes},
   which is especially important
   for polymorphisms.

\subsection{Markov operators}

   The functional analog of the notion of  polymorphism is the notion of
   {\it Markov operator}
   in some functional space, which in this paper will be the Hilbert space
   $L^2_{\mu}(X)$.

\begin{definition}
  {A linear operator $V$ in $L^2_{\mu}(X)$ is called a Markov operator
   if

  {\rm1)} $V$ is a contraction: $\|V\|\leq 1$ in the operator norm;

  {\rm 2)} $V{\it 1}= {V}^*{\it 1}={\it 1}$;\footnote{This condition
  expresses the invariance of the measure
    under the polymorphism (see below); the equality $V^*{\it 1}={\it 1}$ follows automatically from
    the other conditions.}

  {\rm 3)} $V$ is positive, which means that
  $Vf$ is a nonnegative function provided that  $f\in L^2_\mu(X)$ is nonnegative.}
\end{definition}

  It is easy to prove that the set $\cal M$
  of all Markov operators is a convex weakly compact semigroup with
  involution $V \to V^*$.

  Unitary (isometric) Markov operators are precisely
  the operators generated by measure-preserving
  auto(endo)morphisms. We generalize this correspondence to polymorphisms.

\begin{statement}
{\rm 1.} Let $\Pi$ be a polymorphism of
a space $(X, \mu)$ with invariant measure;
 then the formula
 $$(W_{\Pi} f)(x)= \int_X f(y)\mu^x(dy)$$
 defines correctly a Markov operator in $L^2$.

 {\rm 2.} Every Markov operator $W$ in the space $L^2_\mu(X)$, where
 $(X,\mu)$ is a Lebesgue space with continuous finite measure,
 can be represented in the form $W=W_{\Pi}$, where $\Pi$ is
 a  polymorphism of $(X,\mu)$ with invariant measure.

{\rm 3.} The correspondence $\Pi \mapsto W_{\Pi}$ is an antiisomorphism
 between the semigroup with involution of $\bmod 0$ classes of
 polymorphisms and the semigroup of Markov operators; this correspondence
 is also an isomorphism of convex compact spaces.
\end{statement}

 The proof follows from the standard theorems of functional analysis (see \cite{KA}),
 and we mention only the formula for the inverse map from the semigroup
 of Markov operators to the semigroup of
 polymorphisms (for more details, see \cite{Vpol}). Let $U$ be a Markov
 operator; the bistochastic measure of the corresponding polymorphism
 is defined as follows:
  $$\mu(B_1 \times B_2)=<U\chi_{B_1},\chi_{B_2}>.$$
 The check of all required assertions is automatic.

 Note that the correspondence $\Pi \mapsto W_{\Pi}$ is a contravariant
 correspondence and reverses the arrows.

 We will denote by $W_{\Pi}$, $W_{\Lambda}$, {\ldots}  the Markov operators
 corresponding to polymorphisms $\Pi$, $\Lambda$, {\ldots}. If $\Pi$
 is an automorphism, then the operator $W_{\Pi}$ is unitary.
It is clear that this correspondence
 extends the ordinary correspondence
 $T\mapsto U_T$, $(U_Tf)(x)=f(T^{-1}x)$,
 between measure-preserving automorphisms
 and unitary multiplicative real operators (= automorphisms of the unitary ring).

 The compact space $\cal M$ is the convex weakly closed hull of the group of
 unitary multiplicative real operators.

 The orthogonal projector \textbf{1} onto the one-dimensional subspace
 of constants is the Markov operator corresponding to the zero polymorphism:
 $\textbf{1}=W_{\Theta}$.

 An equivalent and more useful definition of the Markov operator corresponding to
 a polymorphism is as follows.
 Consider a bistochastic measure $\nu$ on the space $X \times X$
 and the Hilbert space $L^2_{\nu}(X\times X)$. Consider two subspaces
 $H_1$ and $H_2$ in this space that are the images of $L^2_{\mu}(X)$ under the embedding
 of the spaces of functions of the first and second arguments, respectively,
 to the whole space $L^2_\nu(X\times X)$:
 $$H_1=L^2_{\mu}(X)\longrightarrow  L^2_{\nu}(X \times X) \longleftarrow L^2_{\mu}(X)=H_2.$$
 Denote the orthogonal projection onto the subspace $H_i$ by $P_i$,
 $i=1,2$; then the above definition of Markov operators coincides
 with the following one.

\begin{statement}
\begin{eqnarray*}
 W_{\Pi}f&=&P_2\cdot P_1f, \qquad f\in H_2;\\
 W_{\Pi^*}g\equiv(W_{\Pi})^*g &=& P_1 \cdot P_2 g,\qquad g\in H_1.
\end{eqnarray*}
 \end{statement}

\smallskip

  It is worth mentioning that the conditional expectation, or
  orthogonal projection, onto the subalgebra of functions that
  are constant on elements of a measurable partition
  is the Markov operator corresponding to the ``tautological''
  polymorphism, which we defined in the previous section.

  \subsection{Properties of polymorphisms in terms of Markov operators}

It is not difficult to reformulate all the notions introduced for polymorphisms (ergodicity,
mixing, primality, density, etc.) in terms of Markov operators.

First of all, the mean and pointwise ergodic theorems for Markov operators have the following form
(this is an old result by
   Hopf and Chacon--Ornstein, see \cite{Nev}):
$$\lim_{n\to \infty}\frac{1}{n}\sum_{k=0}^{n-1}(W^k_{\Pi}f)(x)=Pf,$$
where $P$ is the projection onto the maximal fixed
subspace (subalgebra); $P=\Theta$ if the polymorphism
$\Pi$ is ergodic.

A Markov operator $W$
corresponds to a {\it mixing polymorphism} $\Pi$
 if and only if the sequence $W^n$ weakly tends, as $n \to \infty$,
 to the projection onto the subspace of constants:
 $$W^n \to \textbf{1}=W_\Theta.$$
 We will discuss this property in detail later.

We will call a Markov operator $W=W_{\Pi}$ {\it dense} if the
$W$-image of $L^2_{\mu}(X)$ is dense in
$L^2_{\mu}(X)$. Obviously, $W=W_{\Pi}$ is dense if and only if the corresponding
polymorphism $W$ is dense.  The density of $W$
is equivalent to the following condition: the conjugate
operator $W^*$ has zero kernel. In the terminology of the book \cite{FN},
a contraction with dense image is
called a quasi-affinitet.

\begin{definition}
{A Markov operator $V$ is called totally nonisometric
if there is no nonzero subspace in the orthogonal complement
to the subspace of constants in
$L^2_{\mu}(X)$ on which $W$ is isometric.}
\end{definition}

\begin{statement}
{A Markov operator is totally nonisometric if
and only if it corresponds to a prime polymorphism.}
\end{statement}

We are interested in
Markov operators that are far
from isometries
(in other words, in polymorphisms that are far from
automorphisms).
Of course, a mixing Markov operator is totally nonisometric,
but at the same time it is not true that every
Markov operator is the direct sum of a
mixing and isometric Markov operators.
Our main examples will illustrate this effect.

Recall the following notation from operator theory (see \cite{FN}),
which we use for the case of Markov operators.

\begin{definition}
   {A contraction $W$ acting in a Hilbert space $H$ belongs to the classes
   $C_{0,\cdot}$, $C_{\cdot ,0}$, $C_{1,\cdot}$, $C_{\cdot,1}$
   if for every function $f \in H$ that is orthogonal to the subspace of
   constants we have $W^n f\to 0$, $W^{*n}f \to 0$, $W^n f \nrightarrow 0$,
   $W^{*n} f\nrightarrow 0$, respectively.
   The classes $C_{a,b}$, $a,b=0,1$, are defined in
   an obvious way.
   All these classes are nonempty.}
\end{definition}

 We are interested mainly in Markov operators in $L^2_\mu(X)$
 of the class $C_{1, \cdot}$ (or  $C_{\cdot,1}$, or
 $C_{1,1}$),  which
 are totally nonisometric and, consequently,  correspond  to nonmixing prime
 polymorphisms (respectively, polymorphisms whose conjugates
 are prime nonmixing; or polymorphisms such that both the polymorphism and its
 conjugate are prime nonmixing). This class is also the most interesting
 from the viewpoint of the pure operator theory of
 contractions.
 In \cite{FN} it was proved that contractions of type  $C_{1,1}$ are quasi-similar to unitary
 operators; we will extend this fact to Markov operators in Section 4.

 The existence of totally nonisometric nonmixing Markov operators
 is not  {\it a priori} obvious.
 We will describe all such examples. The main feature of such examples is that they are
 not the direct products of mixing and pure deterministic operators.

 The convex structure of the set of all polymorphisms $\cal P$ and
 the isomorphic compact set
 of Markov operators $\cal M$ is very important. Isometries and unitary
 operators are extreme
 points of this compact set, but there are many other extreme points.
 In \cite{Vpol} it was proved that the set
 of extreme polymorphisms is an everywhere dense $G_{\delta}$-set in $\cal P$.
 These extreme Markov
 operators (polymorphisms) have many interesting properties (see
 \cite{Vpol}, \cite{Sud}).

  From the viewpoint of the theory of $C^*$-algebras, it is natural to consider
 the $C^*$-algebras generated by some class of multiplicators (say,
 continuous functions)
 and a given Markov operator and its congugate; this is a
 generalization of the ordinary notion of cross product
 (with an action of the group $\Bbb Z$) and
 cross products with endomorphisms (see a recent paper \cite{ExVe}).
 One of the open questions
 concerns the amenability of the corresponding  $C^*$-algebra.

  Let us also mention the following problem, which was formulated in \cite{Versib}:

\medskip\noindent
 \textbf{Problem.}
  To characterize the $C^*$-algebra $\rm{Alg}({\cal M})$
  generated by all Markov operators in
  $L^2_{\mu}(X)$. This (nonseparable) algebra does not coincide with
  the algebra $B(L^2_\mu(X))$ of all bounded
  operators.\footnote{For example, in $L^2(\Bbb R, m)$ (where $m$ is
  the infinite Lebesgue measure), the operator of Fourier
  transform does not belong to this algebra, as was observed by
  G.~Lozanovsky, see \cite{Versib}.}
  On the other hand,
  this algebra is distinguished and plays the same role in measure
  theory and the theory of Markov operators as
  the algebra of bounded operators $B(H)$ plays in operator
  theory.

   \subsection{The operator formulation of  quasi-si\-mi\-la\-ri\-ty}

\begin{definition}
{A Markov operator $W_1$ is called a quasi-image of a  Markov operator
$W_2$ if there exists a dense Markov operator $U$ such that $U\cdot
W_1=W_2\cdot U$. Markov operators $W_1$ and $W_2$ are called quasi-similar
if each of them is a quasi-image of the other.}
\end{definition}

As follows from definitions, two polymorphisms $\Pi_1, \Pi_2$
are quasi-similar if and only if the corresponding Markov operators
are quasi-similar. The same is true for quasi-images.

   Now let us define more accurately
   the problem of {\it quasi-similarity of automorphisms and polymorphisms},
   which was formulated above in terms of operator formalism.
   Denote by $U_T$ the unitary operator in $L_{\mu}^2(X)$ corresponding
   to an automorphism $T$ and by $W_{\Pi}$ the Markov operator in $L_{\mu}^2(X)$ corresponding
   to a polymorphism $\Pi$; as we have seen, $\Pi$ is prime if and only if
   $W_\Pi$ is {\it totally nonisometric}, i.e., $W_\Pi$ has no invariant subspaces
   (except the one-dimensional subspace of constants) on which it acts as an isometry.
   We also may assume without lost of generality that
   $W_{\Pi^n} \nrightarrow P$ (in the weak topology; here $P$ is the projection
   to the subspace of constants), because in our case $\Pi$ is nonmixing.
   Under this condition, our problem is formulated as follows:

   \smallskip
   \textbf{Problem}. When the Markov operator $U_T$ can be a quasi-image of
   the Markov
   operator $W_\Pi$ and vice versa? When they are quasi-similar
   in the sense of the previous definition?

\section{Markov processes associated with polymorphisms and dilations
of Markov operators}

\subsection{Markov processes}
   Let $\Pi$ be a polymorphism of a space $(X,\mu)$ with invariant measure,
   and let $\nu$ be the corresponding bistochastic measure on
   $X\times X$. As we have mentioned above, every polymorphism with
   invariant measure generates a stationary Markov process;
   thus we consider $\nu$ as the {\it two-dimensional distribution}
   of a stationary Markov process $\Xi(\Pi)$. For this process,
   $X$ is the state space and $\mu$ is an invariant one-dimensional
   distribution. Denote by $M$ the Markov measure in the space
   ${\cal Y} \equiv X^{\Bbb Z}$ of realizations
   of the process  $\Xi(\Pi) \equiv \{\xi_n,\,n \in \Bbb Z\}$,
   and by $S=S_{\Pi}$ the right shift in the measure space
   $({\cal Y}, M)$. The space $L^2_{\mu}(X)$ is naturally embedded
   into the space $L_{M}^2(\cal Y)$ as the subspace of functionals
   of realizations of the process
   that depend only on the value $\xi_0$ of the process at time zero.
   The state space at time $n$,  which is identical to $X$,
   will be denoted by $X_n$.
   Let $S$ be the right shift in the space $({\cal Y},M)$;
   it preserves the measure $M$ and is called the {\it Markov shift}
   corresponding to the polymorphism $\Pi$. For example, if $\Pi=\Theta$,
   then the Markov shift is a Bernoulli shift.

   Note that the Markov process corresponding to
   the conjugate polymorphism $\Pi^*$ is obtained from
   the Markov process of $\Pi$ by reversing time.

   Recall (see \cite{FN}) that if $W$ is a contraction acting in a subspace
   $L \subset H$ of a Hilbert space $H$, then a unitary operator $U$ acting in $H$
   is called a {\it dilation} of $W$ if $$ W^n=PU^n,$$ where $P$ is
   the orthogonal projection $P:L \to H$.
   Every contraction has the so-called minimal dilation (see \cite{FN}).

\begin{definition}
   A dilation of a Markov operator $W$ in
   $L^2_{\mu}(X)$ is a Markov operator $U$
   in some space $L^2_{\alpha}(Y)\supset L^2_{\mu}(X)$
   such that $$ W^n=PU^n,$$ where $P$ is
   the (positive) orthogonal projection $P:L \to H$.
\end{definition}

\begin{statement}
   {The unitary operator $U_S$ in $L_{M}^2(\cal Y)$
   is a Markov dilation of the Markov
   operator $W_{\Pi}$, which acts in the space $L^2_\mu(X_0)$
   regarded as a subspace of
   $L_{M}^2(Y)$:
   $$W_{\Pi}=P_0 U_S.$$
   Here $P_0$ is the expectation (orthogonal projection) onto $L^2_\mu(X_0)$.
   This dilation is not the minimal dilation of $W_{\Pi}$
   in the sense of operator theory, but it is the minimal \textbf{Markov}
   dilation (see \cite{Vpol}).}
\end{statement}

   A general problem is to characterize invariant properties of the Markov shift $S$
   in terms of the polymorphism (= Markov generators), for example, to give a
   characterization of the Bernoulli and non-Bernoulli properties of Markov shifts,
   or to describe relations between regular Markov processes and $K$-automorphisms, etc.

\subsection{Mixing, primality, and tail $\sigma$-field of Markov processes}

   It is clear that the ergodicity of a polymorphism $\Pi$ is equivalent
   to the ergodicity of the process $\Xi(\Pi)$ and to the ergodicity
   of the Markov shift $S$, which is an invariant property.

   Contrary to this,  mixing and primality and other properties of
   a polymorphism and the
   corresponding Markov process are not invariant
   properties of the Markov shift regarded
   as an abstract measure-preserving transformation,
   but can vary for different generators.

   Assume that $\Pi$ is a polymorphism of $(X,\mu)$ and
   $\Xi=\Xi(\Pi)=\{\xi_n\}_{n \in \Bbb Z}$ is the corresponding Markov process
   with state space $X$ and Markov measure $M$ in
   ${\cal Y} =X^{\Bbb Z}$. Denote by  ${\mathfrak{A}}_n$ the
   $\sigma$-subfield in $\cal Y$
   generated  by the set
   of one-dimensional cylindric sets at time $n$,
   and by ${\mathfrak{A}}_-$ (respectively, ${\mathfrak{A}}_+$)
   the tail $\sigma$-field of the past
   (respectively, future) of the process.\footnote{Or the
   intersection over all positive $n$
   of the $\sigma$-fields generated by
   the values of the process before time $-n$
   (respectively, after time $n$).}
     Denote the corresponding partitions
   into infinite pasts (futures)
   by $\xi_{\mp \infty}$, and the quotient spaces
   with measures (the ``infinite past'' and the ``infinite future'')
   by $(X_{\mp},M_{\mp})=(X^{\Bbb Z},M)/{\xi_{\mp \infty}}$.
   These spaces can also be called the
   {\it infinite entrance boundary} and {\it exit boundary}.

   Let $S=S_{\Pi}$ and $S^{-1}$ be the right and left shifts, respectively,
   in the space $({\cal Y}, M)$.

    Recall that a stationary process (even not necessarily Markov)
    is called {\it regular} or {\it pure nondeterministic} in the past (future)
    if the tail $\sigma$-field of the past (future)  (or the entrance (exit) boundary)
    is trivial.\footnote{The term ``regularity''  in this
    sense was first used by Kolmogorov.} Let us emphasize that this is not an
    invariant property of the shift, but a property of the generator (process).

    If the tail $\sigma$-field
    is not trivial, then almost every point $x_\mp\in X_\mp$ determines
    the conditional Markov (nonhomogeneous in time) process $\{\xi_n^{x_\mp}\}$.
    The correspondence $x_\mp\mapsto\{\xi_n^{x_\mp}\}$ determines the
    decomposition of the whole space $(\cal Y, M)$ and the process $\{\xi_n\}$ into a
    direct integral over the quotient space $X_{\mp}$.
    We can correctly define the conditional measure on the
    $\sigma$-field ${\mathfrak{A}}_0$ as the one-dimensional
    distribution at moment $0$ of the conditional Markov process.

    Now we can summarize the information on the tail $\sigma$-fields
    of Markov processes and mixing in the following theorem. We formulate it
    only for the $\sigma$-field of the past; the same is true for the future.

\begin{theorem}
   {For a Markov process
    $\Xi=\Xi(\Pi)$, the following assertions are equivalent:

    1) The process $\Xi$ is regular
    in the past, which means that the $\sigma$-field
    ${\mathfrak{A}}_{-}$ is trivial, i.e.,
    $X_-$  is a one-point space.

    2) The limit (which exists with probability 1) of the conditional measures
    $$\lim_{n\to -\infty} \Pr\{a\mid x_n\}, \quad
    a\in{\mathfrak{A}}_0, \quad x_n \in X_n,
    $$
    on  $\mathfrak{A}_0$
    does not depend $\bmod 0$ on the trajectory $\{x_n\}_{n \in \Bbb Z}$
    and coincides with the unconditional measure.

    3) The Markov generator $\xi_0$ is a $K$-generator
    for the right shift $S$;

    4) The polymorphism $\Pi$ is mixing, i.e., $\Pi^n \to \Theta $.}
\end{theorem}

    The equivalence of the first three claims follows
    more or less from definitions
    and, in contrast to the equivalence with claim 4),
    does not use the Markov property.
    The equivalence between 3) and 4) for Markov processes
    is well known and can be proved directly.

    For Markov chains,  i.e., processes with {\it finite
    or countable state space}, and for some special cases of polymorphisms
    (= transition probabilities), the mixing property is equivalent to
    the property that is usually called the ``absence of  nontrivial
    subclasses of the state space'' (or irreducibility, or convergence
    of the powers of the transition matrix to an invariant vector, etc.;
    see, e.g., {\rm\cite{Kai}}) and that in this paper we have called
    ``primality'' --- the absence of nontrivial
    factor endomorphisms --- or, equivalently, to the following property: a
    Markov process
    has no nontrivial deterministic quotients.

    Thus for Markov chains with discrete state space, we can add to the
    above theorem the following fifth claim, which is equivalent to 1)--4)
    in this case, but is not equivalent to them in the general case:

 \smallskip
    5) {\it There are no nontrivial partitions of the state space
    invariant with respect to the matrix of
    transition probabilities} (see \cite{Vin}).

 \begin{definition}
{A markov process corresponding to a prime nonmixing polymorphism will be called
 quasi-deterministic.}
 \end{definition}

\smallskip
    M.~Rosenblatt was perhaps the first to point out
    the existence of quasi-deterministic (``paradoxical'') Markov processes
    (see \cite[4.4]{Ros}). A more general construction for Anosov systems was
    suggested in \cite{Vgrub} and was called {\it superstability}.

  Thus such a process is not regular, but has no deterministic factors.
  Note that in the case of a quasi-deterministic Markov process
  any measurable set from the $\sigma$-fields
  $\mathfrak A_{\mp}$ of measure not equal to $0$ or $1$ is not a cylindric set (see also \cite{Vgrub}).
  This also contrasts with the theory of Markov chains, where such
  sets are one-dimensional cylinders.
  In the last section we give a description of quasi-deterministic Markov
  processes; it turned out that precisely these polymorphisms are quasi-similar
  to $K$-automorphisms.

\section{Quasi-similarity of automorphisms and polymorphisms and $K$-property}

\subsection{The structure of a quasi-deterministic Markov process}

   Here we briefly describe the structure of the past of a quasi-deterministic
   Markov process.

   It is well known (see, e.g., \cite{Rokh}) that for any stationary process
   with discrete time (even not necessarily Markov)  there is a {\it canonical
   automorphism that acts on the tail $\sigma$-fields ${\mathfrak{A}}_{\pm}$
   and on the quotient spaces $X_{\pm}$}:
   this is the restriction of the left shift $S^{-1}$
   (respectively, right shift $S_+$) to this $\sigma$-field
   and, consequently, to the quotient spaces $X_{\pm}$; it is called the
   {\it residual (tail) automorphism}.\footnote{In ergodic theory, it is
   sometimes called the Pinsker automorphism, and the
   $\sigma$-field ${\mathfrak{A}}_{\pm}$ is called the Pinsker $\sigma$-field.}
   
   Denote these automorphisms by $S_{\mp}$ (it is convenient to use the
   left shift $S^{-1}$ in the past and the right shift $S$ in the future).
   Now we can formulate the first theorem on interrelations between
   the past and present.

\begin{theorem}
   {Assume that $\Pi$ is a prime nonmixing polymorphism.
   Let  $\Xi=\{\xi_n\}_{n \in \Bbb Z}$ be the quasi-deterministic
   stationary Markov process associated with $\Pi$. Then

{\rm 1.} The tail $\sigma$-field ${\mathfrak{A}}_{-}$ of the process is not
   trivial. The tail (residual) automorphism $S_{-}$ acting on the space
   $(X_{-},M_-)$ is ergodic.

{\rm2.} Define a polymorphism $\Lambda'$ as follows:
   the value $\Lambda(x_{-})$, $x \in X_{-}$,
   is the conditional measure $\mu_{x_{-}}$ on the state space $X_0$
   corresponding to the point $x_{-}$ of the tail space $X_{-}$.
   Then
    \begin{equation}
    \Pi\cdot \Lambda' = \Lambda'\cdot S_{-},
\label{1}
    \end{equation}
    in other words, $S_{-}$ is quasi-similar to the polymorphism $\Pi$.

    The polymorphism $\Lambda'$ is injective (see the definition in
    Section 2).

    {\rm3.} The conjugate polymorphism $\Lambda'^*$ from the space $(X_0, \mu)$ to
    the tail space  $(X_{-}, M_{-})$ is also injective; its value at a point $x \in X_0$
    is the conditional measure  on the infinite past $X_{-}$ given that
    the value of the process
    at zero time is equal to $x$.

    {\rm4.} There exists an isomorphism $Q$ between the state space $(X_0,\mu)$
   and the infinite past (tail space) $(X_{-},M_{-})$,
   which determines an automorphism $T$ of the state space $(X_0,\mu)$
   by the formula $T=Q^{-1} S_{-}Q$ and a polymorphism of this space
   by the formula $\Lambda= \Lambda' Q$ such that the automorphism $T$
   is quasi-image of the polymorphism $\Pi$
     :
$$\Lambda\cdot T=\Pi\cdot \Lambda.$$}
\end{theorem}

    In our terminology, the last formula means that $T$ is a quasi-image of
    $\Pi$; if the polymorphism $\Pi^*$ is also prime and nonmixing, then we have
    an analogous formulation with the tail $\sigma$-field of the ``future''
    and obtain the
    {\it quasi-similarity} between $T$ and $\Pi$.
    We omit the proof and observe that the main part of the theorem
    is the ``Markov'' or ``ergodic'' analog of the
    corresponding theorem on contractions
    in Hilbert spaces (see \cite[Ch.~2]{FN}) with some serious complications.
    Indeed, we consider two subspaces in $L^2(\Xi,M)$:
    the space of functions measurable
    with respect to the tail $\sigma$-field and the state space at zero moment;
    these subspaces generate the polymorphism $\Lambda'$ (see Section 2).

\begin{corollary}
   {We can express the polymorphism $\Lambda$ directly in terms of the main
   ingredients $T$ and $\Pi$:
\begin{equation}
\Lambda=\lim_{n\to \infty}\Pi^n T^{-n}.
\label{3}
\end{equation}}
\end{corollary}

   The last theorem reduces the quasi-similarity between an
   automorphism and a polymorphism to the
   state space $(X,\mu)$; the role of the tools of the theory of
   Markov processes is simply in using the
   residual automorphism and interlacing polymorphisms (conditional measures).
   Below we will prove that $T$ is a $K$-automorphism.

   For our purposes, it
   is very convenient to represent
   the polymorphism $\Pi$ as the product $$\Pi=\Phi\cdot T,$$
   where $\Phi$ is a new polymorphism;
   then it is easy to check that
   $$\Lambda=\lim_{n\to \infty}\Phi\cdot T\Phi T^{-1}\ldots
   T^n \Phi T^{-n};$$

   If we set
       $\Phi_k=T^k\Phi T^{-k}$, $k\in \Bbb Z$, $\Phi \equiv \Phi_0$, then
   $$\Lambda=\lim_{n\to \infty} \prod_0^n \Phi_k=\prod_0^{\infty}\Phi_k. $$

   In this setting, the main formula $\Pi \cdot\Lambda = \Lambda\cdot T$
   takes the following form:
   $$\Phi_0 \cdot T = \prod_0^n \Phi_k = \prod_0^n \Phi_k \cdot T.$$

   The convergence of the infinite product is the only condition for
   the construction to be well-defined; in the above situation, this
   follows from
   the existence of the conditional measures on $X_0$ with respect to the infinite past.
   In the examples of the next section it will be proved directly.

   As a result of this section,
   {\it for every nonmixing prime polymorphism $\Pi$ of a space
   $(X,\mu)$ we have found (using the corresponding
   Markov process) an automorphism $T$  such that
   $$\lim_{n\to \infty}\Pi^n T^{-n}=\Lambda$$ and equation
   \rm{(\ref{1})} holds}. The Markov property of the processes was used
   to reduce the problems to a single space (the state space).

\smallskip\noindent
   {\bf Remark.} The orbit partition of the polymorphism $\Lambda$
      is an analog of the partition into subclasses in the theory
   of Markov chains. From the above formulas we can conclude that
   {\it it is an ergodic (absolutely nonmeasurable) equivalence relation}
   (the same is true for the polymorphism $\Lambda^*$).
   Indeed, on the one hand, its measurable hull is an invariant partition
   for $\Pi$, hence it is the trivial partition; but on the other hand, the
   orbit partition is not equal to the trivial partition,
   because the polymorphism $\Pi$ is nonmixing;
   consequently, this is an absolutely nonmeasurable partition.
   We will use this property in the next section.
   This fact is crucial; it illustrates our previous remark on the role of
   the nonmeasurability of the partition into subclasses
   for general Markov processes.

  The Markov shift we have considered obviously has a positive entropy and
  under natural assumptions is a $K$-automorphism (but the Markov
  generator is not a $K$-generator!).
  We will discuss this elsewhere.

\subsection{Random perturbations of $K$-automorphisms}

Using the symbolic representation of $K$-automorphisms (see, e.g., \cite{KSinF}),
we can give a generalization of the construction from \cite{Vgrub},
which associates with
any $K$-automorphism $T$ a polymorphism that is a quasi-image of $T$,
and also find polymorphisms for which $T$ is a quasi-image.
For the special case of Bernoulli automorphisms, it is possible to find
polymorphisms that have both properties simultaneously and, consequently,
are quasi-similar to $T$.

 \begin{theorem}
   {For every $K$-automorphism $T$ there exists a prime nonmixing polymorphism
   $\Pi_{-}$ such that $T$ is a quasi-image of $\Pi_-$, and a polymorphism $\Pi_{+}$
   that is a quasi-image of $T$. If there exists a symbolic realization of $T$
   with finite or countable state space in which the homoclinic equivalence
   relation is ergodic,\footnote{Recall that two sequences
    $\{x_i\}$ and $\{y_i\}$ of the space of sequences
 with a shift-invariant  measure belong to the same
 homoclinic class if for sufficiently large $N$ we have $x_i=y_i$ for $i>N$
 (see \cite{Gord,Vgrub}). It is not known to the author whether
 the homoclinic equivalence relation is ergodic for an arbitrary
  $K$-automorphism.}
  then there exists a polymorphism that is quasi-similar
   to $T$; this is the case of Bernoulli automorphisms.}
\end{theorem}

    The construction is more or less direct, the proof includes some
    combinatorial construction;
    For a given $K$-automorphism $T$ the set of such polymorphisms
    is very large but it is possible to
    describe all such polymorphisms in terms of symbolic realization of $T$
    - these are is nothing more than a random perturbations of the
    $K$-automorphism, along  ``stable,''  ``unstable,''
    or homoclinic
    (which is the intersection of a stable and an unstable one,
    if it is nonempty) leaves of the
    automorphism, respectively;
    each case gives the corresponding quasi-images or quasi-similarity.
    In \cite{Vgrub}, the following terminology was used: a
    $K$-automorphism is ``super-stable'' in the past (respectively,
    in the future, or both), the reason for the term is the following: after
    a those random perturbations, nevertheless the initial automorphism
    can be recovered up to isomorphism as the tail automorphism in the past or future.
    In short, this means that {\it the random perturbation
    of a $K$-automorphism allows us to recover automatically the original automorphism.}

\section*{Acknowledgments}
\addcontentsline{toc}{section}{Acknowledgments}

During a long time I have had a lot of useful discussions on the
topic considered in this paper with my colleagues, especially with
D.~Arov, B.~Rubstein, M.~Rosenblatt, L.~Gandelsman,
V.~Kaimanovich, L.~Khal\-fin, and M.~Gordin; to all of them I am
grateful for information on the literature and related subjects. I
am very grateful to N.~Tsilevich for help with preparing this
text.


\begin{thebibliography}{18}

\addcontentsline{toc}{section}{Bibliography}

\bibitem{Kai} K.~L.~Chung, {\it Markov Chains with Stationary Transition
Probabilities}.
Springer-Verlag, Berlin--G\"ottingen--Heidelberg, 1960.

\bibitem{ExVe} R.~Exel and A.~Vershik, C*-algebras of irreversible dynamical
systems. {\it Canad. Math. J.}, to appear.

\bibitem{Mis} S.~Goldstein, B.~Misra, and M.~Courbage,
On intrinsic randomness of dynamical systems.
{\it J. Stat. Phys.} \textbf{25} (1981), No. 1, 111--126.

\bibitem{Gust} R.~Goodrich, K.~Gustafson, and B.~Misra, On $K$-flows
and irreversibility. {\it J. Stat. Phys.} \textbf{43} (1986), No. 1/2,
317--320.

\bibitem{Gord} M.~I.~Gordin, Double extensions of dynamical systems and
the construction of mixing filtrations, {\it Zap. Nauchn. Semin.
POMI} {\bf 244} (1997), 61--72. English translation: {\it
J. Math. Sci., New York} {\bf 96} (1999), No. 5, 3493--3495.

\bibitem{KA} L.~V.~Kantorovich and G.~P.~Akilov, {\it Functional
Analysis}. Nauka, Moscow, 1977.

\bibitem{KSinF} I.~B.~Kornfeld,  Ya.~G.~Sinai, and S.~V.~Fomin, {\it
Ergodic Theory}. Nauka, Moscow, 1985.

\bibitem{Lax} P.~Lax and R.~Phillips, {\it Scattering Theory}.
Academic Press, New York--London, 1967.

\bibitem{Pr} B.~Misra and I.~Prigogine,
Time, probability, and dynamics. In:
{\it Long Time Prediction in
Dynamics}, C.~W.~Horton, Jr. et. al. (eds.).
Wiley, New York, 1983, pp. 21--43.

\bibitem{FN} B.~Nagy and C.~Foias,
{\it Harmonic Analysis of Operators on Hilbert Spaces}.
Akad\'emiai Kiad\'o, Budapest; North-Holland Publishing Company,
 Amsterdam--London, 1970.

\bibitem{Nev} J.~Neveu, {\it Mathematical Foundations
of the Calculus of
Probability}. Holden-Day, San Francisco, 1965.

\bibitem{Ren} J.~Renault, {\it A Groupoid Approach to $C*$-Algebras}.
Lecture Notes
in Math., vol. 793. Springer-Verlag, Berlin--Heidelberg--New York, 1980.

\bibitem{Rokh} V.~A.~Rokhlin, Lectures on the entropy theory
of measure-preserving transformations. {\it Usp. Mat. Nauk} {\bf22} (1967),
No. 5(137), 3--56. English translation: {\it Russian Math. Surveys}  {\bf22} (1967), No. 5, 1--52.

\bibitem{Ros} M.~Rosenblatt, {\it Markov Processes. Structure and
Asymptotic Behavior}. Springer-Verlag, Berlin--Heidelberg--New York, 1971.

\bibitem{Sud} V.~N.~Sudakov, {\it Geometric Problems in the Theory of
Infinite-dimensional Probability Distributions}.
Tr. Mat. Inst. Steklov {\bf141} (1976). English translation:
Proc. Steklov Inst. Math. {\bf141} (1979).

\bibitem{Vpol} A.~M.~Vershik,
Multivalued mappings with invariant measure (polymorphisms) and Markov operators.
{\it Zap. Nauchn. Semin. LOMI} (1977) {\bf72}, 26--61.
English translation: {\it J. Sov. Math.} {\bf 23} (1983), 2243--2266.

\bibitem{Vgrub} A.~M.~Vershik,
 Superstability of hyperbolic automorphisms and unitary dilatations of
 Markov operators.
 {\it Vestn. Leningr. Univ., Ser. I}, No. 3 (1987), 28--33.
 English translation: {\it Vestnik Leningrad Univ. Math.} {\bf20} (1987),
 No. 3, 22--29.

\bibitem{Versib} A.~M.~Vershik, Measurable realizations of automorphism
groups and integral representations of positive operators.
{\it Sib. Mat. Zh.} {\bf 28} (1987), No. 1(161), 52--60. English translation:
{\it Sib. Math. J.} {\bf 28} (1987), 36--43.

\bibitem{Valg} A.~M.~Vershik, Theory of decreasing sequences of
measurable partitions. {\it
Algebra i Analiz} {\bf6} (1994), No. 4, 1--68. English translation:
{\it St.~Petersburg Math. J.}  {\bf6} (1995), No. 4, 705--761.

\bibitem{Vin} V.~M.~Vinokurov, Conditions for the regularity of
the stochastic processes. {\it Dokl. Akad. Nauk SSSR} {\bf 113} (1957), No. 5,
950--961.

\end{thebibliography}
\end{document}